\begin{document}
\baselineskip=20pt
\textheight=21cm
\textwidth=17cm
\begin{center}
BEREZIN'S QUANTIZATION ON FLAG MANIFOLDS \\
AND SPHERICAL MODULES\\
A.V.Karabegov\\
Joint Institute for Nuclear Research, Dubna, Russia
\end{center}
{\bf Introduction}

The Berezin's quantization on a symplectic manifold $\Omega$
(see \cite{Ber1}) is given by a family of associative algebras
$\{{\cal A}_h\}$, whose elements are smooth functions on $\Omega$,
together with their representations in Hilbert spaces $H_h$,
where the parameter $h$ plays a role of the Planck constant.
 The multiplication $\ast_h$ in ${\cal A}_h$  must
satisfy  the so-called correspondence principle, as $h\to 0$.

Representation of the algebra ${\cal A}_h$ in $H_h$
maps a function
$f\in {\cal A}_h$ to the operator $\hat f$ in $H_h$. In that case the
function $f$ is called a symbol of the operator $\hat f$.
The concrete examples of quantizations are usually based
on explicit constructions of symbols (see \cite{Ber1}, \cite{Ber2}).

In this paper we show how the theory of spherical Harish-Chandra
modules naturally provides the algebras of covariant and contravariant
symbols on generalized flag manifolds. Morover, it also gives rise to
the algebras of symbols, which so far have no general construction
or definition, as do co- and contravariant symbols.
For all these algebras of symbols we give a general proof of the
correspondence principle. We also prove the conjecture from \cite{CGR}
on rational dependence of the product $\ast_h$ on $h$
on an arbitrary generalized flag manifold.

{\it Acknowledgement.} I wish to express my gratitude to Professors
M.S.Nara\-simhan
and D.Vogan for helpful discussions. I also thank the International
Centre for Theoretical Physics in Trieste, Italy, for their hospitality.

\newpage

{\bf1. The definition of quantization}

We will give now a formal definition of quantization on a symplectic
manifold $\Omega$. Let $F$ be a set of positive numbers with a limit
point 0. For each $h\in F$ let ${\cal A}_h$ be an algebra, which
is a linear subspace of $C^\infty(\Omega)$, with a multiplication
$\ast_h$ and given representation in a Hilbert space $H_h$.
Assume that  ${\cal A}=\cup_h {\cal A}_h$
is dense in $C^\infty(\Omega)$. These data determine a quantization
on $\Omega$, if for $f,g\in{\cal A}$ and $h\to 0$ holds the correspondence
principle,
\begin{equation}
f\ast_h g\to fg,\hskip 0.5cm h^{-1}(f\ast_h g-g\ast_h f)\to i\{f,g\}_\Omega.
\end{equation}
Here $\{\cdot,\cdot\}_\Omega$ is a Poisson bracket on $\Omega$.
The limits in (1) make sense if for each $f\in{\cal A}$
one has $f\in{\cal A}_h$ for sufficiently small values of $h$
(a weak nesting property). This condition is valid if for
$h<h'$ holds ${\cal A}_h\supset {\cal A}_{h'}$  (a strong nesting property).

{\bf 2. Covariant and contravariant symbols}

Suppose a compact Lie group $K$ has an irreducible unitary representation
in a $q$-dimensional vector space $E$ with a hermitian scalar product
$<\cdot,\cdot>$,  $v\in E$ is a unit vector and $dk$ is a Haar measure
on $K$ such that the total measure of $K$ is 1.
For $k\in K$ denote $v_k=kv$. Then the orthogonality relations for
matrix coefficients immediately imply, that for all
$\xi,\eta\in E$ holds the Parseval identity
$$
 <\xi,\eta>=q\int <\xi,v_k><v_k,\eta>dk.
$$
Thus the $K$-orbit of $v$, $\{v_k\}$,
forms a supercomplete system of vectors in the sense of \cite{Ber2},
and one can associate with it the constructions of co- and contravariant
symbols on $K$.  Let $P_k$ be the orthogonal projection operator in
$E$ to the vector $v_k$.

{\bf Definition}. {\it A function $f(k)$ on $K$ is called a covariant
symbol of an operator $A\in End\ E$,
associated with the supercomplete system $\{v_k\}$, if
 $f(k)=tr\ AP_k=<Av_k,v_k>$.}

Such symbols for a compact semisimple Lie group $K$ were studied
in \cite{Wild}.

{\bf Definition}. {\it  A measurable function $g(k)$ on $K$
is called a contravariant symbol of an operator $B\in End\ E$,
associated with the supercomplete system $\{v_k\}$, if
$B=q\int g(k)P_kdk$.}

The proofs of the following two lemmas are trivial.

{\bf Lemma 1.} {\it A measurable function $g(k)$ is a contravariant
symbol of an operator $B\in End\ E$ iff for any operator
$A\in End\ E$ with the covariant symbol $f(k)$ holds
$q\int f(k)g(k)dk=tr\ AB$}.

Let $k_0\in K$. Denote  $v'=k_0v$ and consider the
supercomplete system  $\{v_k'=kv'\}$.

{\bf Lemma 2.} {\it Let $f(k)$ be the covariant (contravariant) symbol
of an operator $A\in End\ E$, associated with the supercomplete system
$\{v_k\}$. Then $f(kk_0)$ is the covariant (contravariant) symbol of
$A$, associated with the supercomplete system  $\{v_k'\}$.}

{\bf 3. A spherical principal series of Harish-Chandra modules}

We give a construction of a spherical principal series of Harish-Chandra
modules, following \cite{Duf}.

Let $\underline{g}$ be a complex semisimple Lie algebra,
$\underline{h}$ its Cartan subalgebra, $W$ the Weyl group,
$\triangle,\ \triangle^+,\ \Sigma$ the sets of all nonzero,
positive and simple roots respectively, $\rho$ the half-sum
of positive roots. For each $\alpha\in\triangle$ choose weight
elements $X_\alpha\in \underline{g}$ such that
$[H_\alpha,X_{\pm\alpha}]= \pm 2X_{\pm\alpha}$ for
$H_\alpha=[X_\alpha,X_{-\alpha}]$.  An element $\lambda\in \underline{h}^*$
is called dominant, if $\lambda(H_\alpha)\geq 0$ for all $\alpha\in\Sigma$,
and is a weight, if $\lambda(H_\alpha)\in{\bf Z}$ for all
$\alpha\in\Sigma$. Denote   $\underline{n},\underline{n}^-$
  the complex subalgebras of $\underline{g}$,
  spanned by $\{X_\alpha\}$ for  $\alpha$ running over
  $\triangle^+,\ -\triangle^+$ respectively. Then
  $\underline{g}=\underline{n}^-\oplus\underline{h}
  \oplus\underline{n}$ is the Gauss decomposition.

  Let $g,h,n,n^-$ denote the Lie algebras
  $\underline{g},\ \underline{h},\ \underline{n}$ and $\underline{n}^-$
  respectively, considered as real Lie algebras,   $\theta$ be an
  automorphism of  $g$, antilinear with respect to the complex structure
  on $\underline{g}$ and such that $\theta X_\alpha=-X_{-\alpha},\ \theta
  H_\alpha=-H_\alpha$. The subspace of $\theta$-invariant elements of
  $g,\ k=g^\theta$ is a compact form of $\underline{g}$. Denote by
  $a$ the real subalgebra of  $\underline{g}$,
  generated by $\{H_\alpha\}$. Then $g=k\oplus a\oplus n$
  is the Iwasawa decomposition.

  Now let $G$ be the complex semisimple connected simply connected
Lie group with the Lie algebra $g$, considered as a real Lie group,
and $G=KAN$  be the corresponding Iwasawa decomposition. Then
 $K$ is a compact semisimple connected simply connected
Lie group. Let $m=k\cap h,\ M$ be a maximal torus in $K$ with the
Lie algebra  $m$, $M'$ its normalizer $K$, so that $W=M'/M$.

  Denote by $\Lambda$  the space of functions on $K$,
$K$-finite with respect to the (left) shifts and invariant
with respect to the right shifts by the elements of $M$
(in \cite{Duf} it is denoted by $L(0)$).
Let $X\in a,\ a=\exp X\in A$. Then for $\lambda\in \underline{h}^*$
$a^\lambda={\rm e}^{<\lambda,X>}$.
Denote by $\tilde\Lambda(\lambda)$ the space of functions $f$
on $G$ such that for $g\in G,m\in M,a\in A,n\in N$ holds
$f(gman)=f(g)a^{2(\lambda-\rho)}$, and the restriction of  $f$ to
$K$ belongs to $\Lambda$ (in the notation of \cite{Duf}, $L(0,-2\lambda)$).
 Define a $(g,K)$-module structure on $\tilde\Lambda(\lambda)$,
 for $f\in\tilde\Lambda(\lambda),\ X\in g,\ k\in K\
Xf(g)=\frac{d}{dt}|_{t=0}f(\exp(-tX)g),\
kf(g)=f(k^{-1}g)$.
     The module $\tilde\Lambda(\lambda)$ is spherical, i.e.,
its subspace of $K$-invariant vectors is one dimensional.

By virtue of the Iwasawa decomposition, the restriction to $K$
determines a bijection of $\tilde\Lambda(\lambda)$ on $\Lambda$.
Denote by $\Lambda(\lambda)$ the space $\Lambda$, endowed with a
$(g,K)$-module structure, induced from $\tilde\Lambda(\lambda)$.
The  $(g,K)$-modules $\Lambda(\lambda)$ (or $\tilde\Lambda(\lambda)$)
form the spherical principal series of Harish-Chandra $(g,K)$-modules.

Let $A(\lambda)$ and $B(\lambda)$ denote respectively the cyclic
submodule of $\Lambda(\lambda)$, generated by $1\in\Lambda(\lambda)$,
and its maximal proper submodule. Then
$\hat\Lambda(\lambda)=A(\lambda)/B(\lambda)$
is the canonical simple subquotient of a principal series.
It is known, that for a given
$\lambda\in \underline{h}^*$
the modules $\hat\Lambda(w\lambda),\ w\in W$, are isomorphic
to each other (see \cite{Duf}).

{\bf 4. The Dixmier mapping}

Let ${\cal U}(\underline{g})$ be the universal enveloping algebra of
$\underline{g}$. Define an action of $g$ on ${\cal U}(\underline{g})$:
for $X\in g,\ u\in{\cal U}(\underline{g})$ let $X:u\mapsto \theta(X)u-uX$,
and let $K$ act on ${\cal U}(\underline{g})$ by the adjoint action
(denoted by $Ad$).
If $X\in k$, then $\theta X=X$, so $X:u\mapsto Xu-uX$, i.e.,
$k$ acts on ${\cal U}(\underline{g})$ by the adjoint
action. Thus the actions of $g$ and $K$ agree on $k$ and
${\cal U}(\underline{g})$ is a $(g,K)$-module.

Consider a morphism of $(g,K)$-modules
$s_\lambda :{\cal U}(\underline{g})\to \Lambda(\lambda+\rho)$,
such that $s_\lambda 1=1$. This requirement defines $s_\lambda$
uniquely, since $1\in{\cal U}(\underline{g})$ is cyclic.
Thus $s_\lambda$ coincides with the morphism
$\Phi_{-\lambda}$, introduced in \cite{Dix}, Sect. 9.6.5.

{\bf Lemma 3.} {\it (i) $Im\ s_\lambda=A(\lambda)$;
         (ii) $Ker\ s_\lambda=Ann\ L(\lambda+\rho)$, where \\
$Ann\ L(\lambda+\rho)$ is the annihilator of the simple
highest weight $\underline{g}$-module \\
$L(\lambda+\rho)$
with the highest weight $\lambda$.}

The proof of (i) is trivial, and (ii) is proved in the lemma 9.6.5,
\cite{Dix}.

{\bf Corollary.} {\it  $A(\lambda)$ carries a natural structure
of an algebra, isomorphic to
${\cal U}(\underline{g})/Ann\ L(\lambda+\rho)$.}

 We will give now an explicit construction of the mapping $s_\lambda$.
  For $\lambda\in \underline{h}^*$ define a functional
$\varphi_\lambda$ on ${\cal U}(\underline{g})$
  (in the notation of \cite{Dix}, Sect.9.6.4., it corresponds to
  $\varphi_{-\lambda}$).
  It follows from the Gauss decomposition, that
    \begin{equation}
   {\cal U}(\underline{g})=(\underline{n}^- {\cal U}(\underline{g})+{\cal
  U}(\underline{g})\underline{n}) \oplus {\cal U}(\underline{h}).
  \end{equation}
  Denote the projection of $u\in {\cal U}(\underline{g})$ on ${\cal
  U}(\underline{h})$ by $u_0$.  Since  ${\cal U}(\underline{h})$
  is naturally isomorphic to the symmetric algebra $S(\underline{h})$
  (or to the algebra of polynomials on $\underline{h}^*$), one can
  evaluate $u_0$ at the point $\lambda\in \underline{h}^*$. Define
  $\varphi_\lambda(u)=u_0(\lambda)$.

     For $u\in {\cal U}(\underline{g})$ define a function
   $\tilde s_\lambda u$ on $K$ by the formula $\tilde s_\lambda
  u(kan)=\varphi_\lambda(Ad(k^{-1})u)a^{2\lambda}$.

{\bf   Proposition 1.} {\it The mapping $\tilde s_\lambda :u\mapsto
  \tilde s_\lambda u$ is a morphism of $(g,K)$-modules from
${\cal U}(\underline{g})$ to $\tilde\Lambda(\lambda+\rho)$.}

  {\it Proof.} The $K$-equivariance of $\tilde s_\lambda$
can be checked immediately. Now the $K$-finiteness of
$\tilde s_\lambda u$ follows from the finiteness of the adjoint action
of $K$ in ${\cal U}(\underline{g})$.  The adjoint action of $M$ on
${\cal U}(\underline{g})$ preserves the decomposition (2)
and is trivial on ${\cal U}(\underline{h})$, which implies
the invariance of $\tilde s_\lambda u$ with respect to
right shifts by the elements of $M$. Thus $\tilde
s_\lambda u\in\tilde\Lambda(\lambda+\rho)$. It remains to show that
$\tilde s_\lambda$ is $g$-equivariant, that is, for
$X\in g,\ u\in{\cal U}(\underline{g})$ the functions $X\tilde s_\lambda
u$ and $\tilde s_\lambda(\theta(X)u-uX)$ coincide. Since the
restriction of $\tilde\Lambda(\lambda+\rho)$ to $K$ is bijective and
$\tilde s_\lambda$ is $K$-equivariant, it is enough to compare these
functions at $e\in K$, that is, to check the equality
 \begin{equation}
 \frac{d}{dt}|_{t=0}\tilde s_\lambda
u(\exp(-tX))= \varphi_\lambda(\theta(X)u-uX)
\end{equation}
for $X\in k,a,n$.  Since the actions of $g$ and $K$ agree on $k$,
(3) holds for $X\in k$.  Let $X\in a$, then $\tilde s_\lambda
u(\exp(-tX))={\rm e}^{-2t<\lambda,X>}\varphi_\lambda(u)$,
so the left-hand side of (3) equals
 $-2<\lambda,X>\varphi_\lambda(u)$.  On the other hand,
 for $X\in a\ \theta X=-X$, so $\theta(X)u-uX=-2uX$.
Since the multiplication by $X\in a$ preserves the
decomposition (2), one has
$\varphi_\lambda(-2uX)=-2<\lambda,X>\varphi_\lambda(u)$. Now let
$X\in n$. Then $\tilde s_\lambda u(\exp(-tX))=\varphi_\lambda(u)$
does not depend on $t$, and thus the left-hand side of (3) is
equal to 0. Since $\theta X\in n^-$, therefore $\theta(X)u-uX$
belongs to the sum in the parenthesis in (2), so the
right-hand side of (3) is also 0. This concludes the proof.

{\bf Corollary.} {\it The mapping $s_\lambda$ is given
by the formula $s_\lambda u(k)=\varphi_\lambda(Ad(k^{-1})u)$.}

{\it Proof.} Since the restriction of $\tilde s_\lambda 1$ to $K$
is identically 1, the morphism
${\cal U}(\underline{g})\ni u\mapsto \tilde s_\lambda u|_K\in
\Lambda(\lambda+\rho)$
coincides with $s_\lambda$. Now the proof follows from the formula of
$\tilde s_\lambda$.

{\bf 5. Symbol algebras}

Let $\lambda\in \underline{h}^*$ be a dominant weight.
Consider a holomorphic irreducible representation of the group
$G$ in $E^\lambda$ with the highest weight $\lambda$, the
$K$-invariant hermitian scalar product $<\cdot,\cdot>$, and the
highest weight vector $v$ of the length 1. Denote the corresponding
representation of ${\cal U}(\underline{g})$ in $E^\lambda$
by $\pi_\lambda$.  Let $w_0\in W$ be the Coxeter element
(the element of the maximal reduced length,
which maps $\triangle^+$ to $-\triangle^+$). Set
$\lambda'=-w_0\lambda$.  It is known, that the module $E^{\lambda'}$
is dual to $E^\lambda$.

Endow the space $End\ E^\lambda$ with a
$(g,K)$-module structure,  for $A\in End\ E^\lambda,$ \\
$u\in{\cal U}(\underline{g})$ and $k\in K$, let
$u:A\mapsto \pi_\lambda(\theta X)A-A\pi_\lambda(X)$ and $k:A\mapsto kAk^{-1}$.

The complexification of the algebra $g$ is isomorphic to
$\underline{g}\times\underline{g}$. Consider an imbedding of $g$ in
$\underline{g}\times\underline{g}$, $X\mapsto (\theta X,X)$. Then
$End\ E^\lambda$ is isomorphic to $E^\lambda\otimes E^{\lambda'}$
as a $\underline{g}\times\underline{g}$-module and therefore is
simple.

It is easy to check that
$\pi_\lambda:{\cal U}(\underline{g}) \to End\ E^\lambda$
is a surjective morphism of $(g,K)$-modules.

{\bf Proposition 2.} {\it Let $\lambda\in\underline{h}^*$ be
 a dominant weight, $u\in{\cal U}(\underline{g})$.
The function $s_\lambda u$ is the covariant symbol of the
operator $\pi_\lambda(u)\in End\ E^\lambda$,
associated with the supercomplete system $\{v_k=kv\}$.}

{\it Proof.} Suppose $Ad(k^{-1})u$ is decomposed according to (2),
  $Ad(k^{-1})u=u_0+Xu_1+u_2Y$, where $u_0\in {\cal
  U}(\underline{h}),\ X\in \underline{n}^-,\ Y\in\underline{n}$.
  Let us calculate the covariant symbol of the operator
  $\pi_\lambda(u)$, associated with $\{v_k\}$:\\

$<\pi_\lambda(u)v_k,v_k>=<\pi_\lambda(Ad(k^{-1})u)v,v>=$
\begin{equation}
<\pi_\lambda(u_0)v,v>+<\pi_\lambda(Xu_1)v,v>+<\pi_\lambda(u_2Y)v,v>.
\end{equation}
  Since $v$ is the highest weight vector of the weight $\lambda$,
the first summand in the right-hand side of (4) is equal to
  $s_\lambda(k)$, and the third one is 0. Finally, since
$\theta X\in \underline{n}$, we have
$<\pi_\lambda(Xu_1)v,v>=-<\pi_\lambda(u_1)v,\pi_\lambda(\theta X)v>=0.$

The proposition is proved.

{\bf Corollary.} {\it Let $\zeta$ be the mapping which maps
an operator $A\in End\ E^\lambda$ to its covariant symbol,
$\zeta A(k)=<Akv,kv>$. Then $\zeta$ is an isomorphism of
$(g,K)$-modules and associative algebras
$End\ E^\lambda$ and $A(\lambda+\rho)$.}

{\it Proof.} Proposition 2 asserts that
$s_\lambda=\zeta\pi_\lambda$. It follows from the surjectivity
of $\pi_\lambda$ and from the fact that the image of $s_\lambda$
coincides with $A(\lambda+\rho)$, that the image of $\zeta$
coincides with $A(\lambda+\rho)$ as well.  The mapping
 $\zeta$ is injective since  $End\ E^\lambda$ is a simple module.
Finally, since $s_\lambda$ and $\pi_\lambda$ are morphisms of
$(g,K)$-modules and homomorphisms of associative algebras,
the same is true for $\zeta$.

 The shifted action of the Weyl group $W$ on $\underline{h}^*$,
  is defined as follows. For $\lambda\in\underline{h}^*,\ w\in W\quad
  w\cdot\lambda=w(\lambda+\rho)-\rho$. For a dominant weight
$\lambda\in\underline{h}^*,\ $ the $(g,K)$-module
  $\hat\Lambda(\lambda+\rho)= A(\lambda+\rho)/B(\lambda+\rho)$,
together with all \\
$\hat\Lambda(w(\lambda+\rho) ),\ w\in W$,
are isomorphic to $End\ E^\lambda$. This implies that, since
the modules  $A(w(\lambda+\rho)),\ w\in W,$ ¨ $End\ E^\lambda$
  are spherical, one can choose a surjective morphism
  $\tau_{w\cdot\lambda}: A(w(\lambda+\rho))\to End\ E^\lambda$
  such that $\tau_{w\cdot\lambda}1=1$. It is clear, that
  $\tau_{w\cdot\lambda}s_{w\cdot\lambda}=\pi_\lambda$,
  and the $(g,K)$-module morphism $\tau_{w\cdot\lambda}$ is also an
  algebra homomorphism, i.e. the representation of
  $A(w(\lambda+\rho))$ in  $E^\lambda$.

  Now we can introduce the following

  {\bf Definition.} {\it A function $f\in A(w(\lambda+\rho))$
  is called a mixed symbol of an operator $A\in End\ E^\lambda$,
  corresponding to the element $w\in W$, if $A=\tau_{w\cdot\lambda}f$.}

  Let $\mu\in\underline{h}^*,\ f\in\Lambda(\mu),\
  g\in\Lambda(-\mu)$, and $q$ be a nonzero constant.  It is known that
  $(f,g)_q=q\int_K f(k)g(k)dk$ is a pairing of $(g,K)$-modules
  $\Lambda(\mu)$ and $\Lambda(-\mu)$ (see \cite{Duf}).
   Since the module $B(\mu)$ does not contain  $1\in A(\mu)$,
  it is orthogonal to $1\in A(-\mu)$, and therefore to the whole
  module $A(-\mu)$, with respect to $(\cdot,\cdot)_q$.
  Besides, $(1,1)_q=q$.  This implies that $(\cdot,\cdot)_q$
  induces a nondegenerate pairing of the simple  $(g,K)$-modules
  $\hat\Lambda(\mu)$ and $\hat\Lambda(-\mu)$, and each pairing
of these modules can be obtained by choosing the appropriate
value of $q$.

  Now let $\lambda\in \underline{h}^*$ be a dominant weight.
  The module $\hat\Lambda(\lambda+\rho)$ is isomorphic to
  $End\ E^\lambda$ and dual to $\hat\Lambda(-\lambda-\rho)$.
  Since  $-\lambda-\rho=w_0(\lambda'+\rho)$, the module
  $\hat\Lambda(-\lambda-\rho)=\hat\Lambda(w_0(\lambda'+\rho))$
  is isomorphic to $End\ E^{\lambda'}$. On the other hand,
  since the $\underline{g}$-modules $E^\lambda$ and $E^{\lambda'}$
  are dual, an operator $A\in End\ E^\lambda$ has the transposed
  $A^t\in End\ E^{\lambda'}$.  It is easy to check, that for
  $A\in End\ E^\lambda$ and $B\in End\ E^{\lambda'},\ trAB^t$
  also defines a nondegenerate pairing of the $(g,K)$-modules
  $End\ E^\lambda$ and $End\ E^{\lambda'}$.  This implies

  {\bf Lemma 4.} {\it Let $f\in A(\lambda+\rho),\
g\in A(w_0(\lambda'+\rho)),$ and $q$
 be the dimension of $E^\lambda$. Then
  $$
  q\int_K f(k)g(k)dk=tr(\tau_\lambda f)(\tau_{w_0\cdot\lambda'}g)^t.
  $$
  }
  To prove the assertion of the lemma it is enough to check it for
  $f=g=1$.

  It follows from Proposition 2, that a function $f$ is the covariant
 symbol of the operator $\tau_\lambda f$, associated with the supercomplete
 system $\{v_k=kv\}$. Now Lemma 1 and  Lemma 4 imply,
 that a function $g$ is a contravariant symbol
of the operator $(\tau_{w_0\cdot\lambda'}g)^t$,
also associated with the supercomplete system $\{v_k\}$.
  Taking in account that for a dominant weight
  $\lambda\quad$ $\tau_{w\cdot\lambda} s_{w\cdot\lambda}=\pi_\lambda$,
  we get the following\\

  {\bf Lemma 5.}  {\it For $u\in{\cal U}(\underline{g})$ the function
  $s_{w_0\cdot\lambda'} u$ is a contravariant symbol of the operator
  $(\pi_{\lambda'}(u))^t\in End\ E^\lambda$,
  associated with the supercomplete system  $\{v_k\}$.}

   Let $u\mapsto\check u$ be an involutive antiautomorphism
  of ${\cal U}(\underline{g})$, such that for
  $X\in\underline{g}\ \check X=-X$. Since the
  $\underline{g}$-modules $E^\lambda$ and $E^{\lambda'}$ are dual,
  the transposed operator of
  $\pi_{\lambda'}(X)\in End\ E^{\lambda'}$ is equal to
  $-\pi_\lambda(X)\in End\ E^\lambda$, which implies that
  $(\pi_{\lambda'}(u))^t=\pi_\lambda(\check u)$.  Therefore
  it follows from Lemma 5, that the function $s_{w_0\cdot\lambda'}u$
  is a contravariant symbol of the operator  $\pi_\lambda(\check u)$,
  associated with $\{v_k\}$.

  {\bf Lemma 6.} {\it Let $\tilde{w}_0\in M'$ be a representative
  of the Coxeter element $w_0\in W$, $u\in{\cal U}(\underline{g})$.
   Then $s_{w_0\cdot\lambda'}u(k)=s_{w_0\cdot\lambda}\check
  u(k\tilde{w}_0^{-1})$.  }

 {\it Proof.}  Since the mappings $u\mapsto\check u$ and
 $u\mapsto s_\lambda u$ are $K$-equivariant,
 it is enough to check the equality of the lemma at the point
 $k=e$.  Let $u$ be decomposed in accordance with (2),
$u=u_0+Xu_1+u_2Y$, where $u_0\in {\cal U}(\underline{h}),\ X\in
\underline{n}^-,\ Y\in\underline{n}$.  Then, taking in account that
$w_0\cdot\lambda'=-\lambda-2\rho$, we get
$s_{w_0\cdot\lambda'}u(e)=u_0(-\lambda-2\rho)$. Since
the antiautomorphism $u\mapsto\check u$ for
$u\in{\cal U}(\underline{h})$ reduces to the change of a sign
in the argument when  ${\cal U}(\underline{h})$ is identified with
$S(\underline{h})$, then
$u_0(-\lambda-2\rho)=\check u_0(\lambda+2\rho)$. On the other hand,
  $\check u=\check u_0-\check u_1 X-Y\check u_2$, so
\begin{equation}
Ad(\tilde w_0)\check u=
Ad(\tilde w_0)\check u_0-Ad(\tilde w_0)
(\check u_1X+Y\check u_2).
\end{equation}
Since $Ad(\tilde w_0)X\in \underline{n}$, and
$Ad(\tilde w_0)Y\in \underline{n}^-$, the subtrahend in the right-hand
side of (5) belongs to the sum in the parenthesis in (2), therefore
  $$
  s_{w_0\cdot\lambda}\check u(\tilde{w}_0^{-1})
=\varphi_{w_0\cdot\lambda}(Ad(\tilde w_0)\check u)=
  (Ad(\tilde w_0)\check u_0)(w_0\cdot\lambda)=
\check u_0(w_0(w_0\cdot\lambda)).
  $$
To conclude the proof it is enough to notice that
  $w_0(w_0\cdot\lambda)=\lambda+2\rho.$

   Consider a vector $v'=\tilde w_0v\in
  E^\lambda$.  Since $v\in E^\lambda$ is the highest weight vector,
  then $v'$ is the lowest weight one.

  {\bf Theorem 1.} {\it (i) A mixed symbol of an operator
  $A\in End\ E^\lambda$, corresponding to the unit element
 $e\in W$, is the covariant symbol of  $A$,
  associated with the supercomplete system  $\{v_k=kv\}$. \\
  (ii) A mixed symbol of an operator  $A$,
       corresponding to the element $w_0\in W$,
  is a contravariant symbol of $A$,
  associated with the supercomplete system $\{v_k'=kv'\}$.}

 {\it Proof.} The assertion (i) is proved in Proposition 2.
  It follows from Lemmas  5 and 6 that for $u\in{\cal U}(\underline{g})$
  the function $s_{w_0\cdot\lambda}u(k\tilde{w}_0^{-1})$
  is a contravariant symbol of the operator $\pi_\lambda(u)$,
   associated with the supercomplete system  $\{v_k\}$. Now it
follows from Lemma 2, that the function  $s_{w_0\cdot\lambda}u(k)$
 is a contravariant symbol of the operator   $\pi_\lambda(u)$,
   associated with the supercomplete system   $\{v_k'\}$,
  which concludes the proof.

{\bf 6. Symbol algebras on flag manifolds}

Denote by $a^*$ a real subspace of $\underline{h}^*$,
formed by those functionals which are real on $a$.
In particular, the roots and weights of $\underline{g}$
belong to $a^*$.
We will give a condition on $\lambda\in a^*$, under which
one can define a pushforward of the algebra $A(\lambda+\rho)$
to an orbit of the adjoint representation of $K$ in $k$.
These orbits are generalized flag manifolds.

The Killing form $(\cdot,\cdot)$ on $\underline{g}$
is positive definite on $a$ and negative definite on $m$.
Let $\lambda\in a^*$. Define a mapping
$\Psi_\lambda:K\mapsto \underline{g}$ such that for each
$k\in K,\ X\in\underline{g}$ holds
$is_\lambda X(k)=(X,\Psi_\lambda (k))$. The following lemma
is straightforward.

{\bf Lemma 7.} {\it For $k,l\in K\quad \Psi_\lambda(k^{-1}l)=
Ad(k)\Psi_\lambda(l)$.}

The lemma asserts, that $\Psi_\lambda$ is a $K$-equivariant mapping.
For each $\lambda\in a^*$ choose $H^\lambda\in a$
such that for $H\in\underline{h}\quad \lambda(H)=(H,H^\lambda)$.
Notice, that for
$\alpha\in\triangle\ H^\alpha=\frac{2}{(H_\alpha,H_\alpha)}H_\alpha$.
Since $ia=m$, for $\lambda\in a^*\ iH^\lambda\in m$.

{\bf Lemma 8.} {\it  Let $\lambda\in a^*$. Then $\Psi_\lambda(e)=iH^\lambda$,
so the image of $\Psi_\lambda$ is the orbit of the point
 $iH^\lambda\in m$ with respect to the adjoint action of $K$ on
 $\underline{g}$, and it lies completely in $k$.}

  {\it Proof.}  For any $X\in\underline{g}$ we have
  $(X,\Psi_\lambda(e))=is_\lambda X(e)=i\varphi_\lambda(X)$.
  If $X\in\underline{n}\oplus\underline{n}^-$, then
  $i\varphi_\lambda(X)=0$ and $(X,H^\lambda)=0$. And if
  $X\in\underline{h}$, then
  $i\varphi_\lambda(X)=i\lambda(X)=(X,iH^\lambda)$.
  Thus  $(X,\Psi_\lambda(e))=(X,iH^\lambda)$ holds for each
  $X\in\underline{g}$,  which implies the assertion of the lemma.

Denote by $\triangle(\lambda)$
the set of roots $\alpha\in\triangle$  such that $\lambda(H_\alpha)=0$
for $\lambda\in a^*$ by $\triangle(\lambda)$, and let
$\triangle^+(\lambda)=\triangle^+\cap\triangle(\lambda)$.
It is evident, that $\lambda(H_\alpha)=0$ iff $\lambda(H^\alpha)=0$.

The centralizer of $H^\lambda$ in $\underline{g}$ is generated
by the Cartan subalgebra $\underline{h}$ and the elements
$X_\alpha,\ \alpha\in\triangle(\lambda)$.
Denote by $k^\lambda\subset k$ the centralizer of $iH^\lambda$
in $k$. Then $k^\lambda$ is generated by the subalgebra $m$ and the elements
$X_\alpha-X_{-\alpha}$ and $i(X_\alpha+X_{-\alpha}),\
\alpha\in\triangle^+(\lambda)$.
Since the stabilizer $K^\lambda\subset K$ of the element
$iH^\lambda\in k$ with respect to the adjoint action of
$K$ is connected, it is the integral subgroup
of the subalgebra $k^\lambda\subset k$.
The orbit $\Omega_\lambda\subset k$ of the point $iH^\lambda$
under the adjoint action of $K$ is isomorphic to $K/K^\lambda$
as a homogeneous space, and is a generalized flag manifold.

Let $\sigma\subset\Sigma$ be a set of simple roots.
Denote by $<\sigma>$  the set of roots which are linear combinations of
elements of $\sigma$.

{\bf Definition.} {\it An element $\lambda\in a^*$
is called regular  relative to $\sigma\subset\Sigma$, if
$\triangle(\lambda)=<\sigma>$.}

{\bf Definition.}  {\it An element $\lambda\in a^*$
is called relatively regular, if it is regular relative to some
$\sigma\subset\Sigma$.}

{\bf Lemma 9.} {\it An element  $\lambda\in a^*$ is relatively
regular iff for all $\alpha,\beta\in\triangle^+$ such that
$\alpha+\beta\in\triangle^+(\lambda)$ holds
$\alpha,\beta\in\triangle^+(\lambda)$.}

{\it Proof.} The necessity is evident. Let us prove the sufficiency.
Set $\sigma=\Sigma\cap\triangle(\lambda)$. We will show that
$\triangle(\lambda)=<\sigma>$. Since
$\triangle(\lambda)=-\triangle(\lambda)$, it is enough to consider
elements $\gamma\in\triangle^+(\lambda)$. If
$\gamma\in\triangle^+$, then either $\gamma\in\Sigma$, or
one can find $\alpha\in\triangle^+,\ \beta\in\Sigma$ such that
$\gamma=\alpha+\beta$. It then follows from
$\gamma\in\triangle^+(\lambda)$ that either $\gamma\in\sigma$,
or $\alpha\in\triangle^+(\lambda)$, and $\beta\in\sigma$.
It remains to use the induction over the number of
summands in the decomposition of $\gamma$
to a sum of simple roots.

{\bf Lemma 10.} {\it Any dominant element $\lambda\in a^*$ is relatively
regular.}

Since the Killing form is positive definite on $a$,  $\lambda\in a^*$
is dominant iff for all $\alpha\in\triangle^+\quad \lambda(H^\alpha)\geq 0$.
Let $\alpha,\beta\in\triangle^+$
be such that $\gamma=\alpha+\beta\in\triangle^+(\lambda)$.
Then $\lambda(H^\alpha)+\lambda(H^\beta)=\lambda(H^\gamma)=0$, and since
$\lambda(H^\alpha)\geq 0,\ \lambda(H^\beta)\geq 0$, then
$\lambda(H^\alpha)=\lambda(H^\beta)=0$, which means that
$\alpha,\beta\in\triangle^+(\lambda)$.
The lemma follows now from Lemma 9.

{\bf Lemma 11.} {\it All the elements of
$A(\lambda+\rho),\ \lambda\in a^*$ are invariant under the right
shifts by the elements of $K^\lambda$ iff
$\lambda$ is relatively regular.}

{\it Proof.} Assume  $\alpha,\beta\in\triangle^+$
are such that $\gamma=\alpha+\beta\in\triangle^+(\lambda)$,
and, say, $\alpha\notin\triangle^+(\lambda)$.
We will show that the function $s_\lambda X_\alpha X_\beta$
is not constant on $K^\lambda\subset K$. In order to do that
let us calculate the infinitesimal action
of the element $X_\gamma-X_{-\gamma}\in k^\lambda$
on $s_\lambda X_\alpha X_\beta$ at the point
$e\in K$:\\

$\frac{d}{dt}|_{t=0}(s_\lambda X_\alpha X_\beta)
(\exp (-t(X_\gamma-X_{-\gamma})))=$
\begin{equation}
\frac{d}{dt}|_{t=0}\varphi_\lambda
(Ad(\exp t(X_\gamma-X_{-\gamma}))X_\alpha
X_\beta)=\varphi_\lambda([X_\gamma-X_{-\gamma},X_\alpha X_\beta]).
 \end{equation}

Since $X_{-\gamma}X_\alpha X_\beta,\ X_\gamma X_\alpha X_\beta$
    and  $X_\alpha X_\beta X_ \gamma$ belong to the kernel of
$\varphi_\lambda$, the expression in (6) is equal to
    $\varphi_\lambda(X_\alpha X_\beta X_{-\gamma})=
    \varphi_\lambda(X_\alpha X_{-\gamma}X_\beta)+
    \varphi_\lambda(X_\alpha[X_\beta,X_{-\gamma}])$.
     Since $\beta-\gamma=-\alpha\in\triangle$, it follows that
    $[X_\beta,X_{-\gamma}]=c\cdot X_{-\alpha}$,
    where $c$ is a nonzero constant. Finally, the expression in (6)
    is equal to $\varphi_\lambda(X_\alpha[X_\beta,X_{-\gamma}])=
    c\varphi_\lambda(X_\alpha X_{-\alpha})=
    c\varphi_\lambda([X_\alpha,X_{-\alpha}])=
    c\varphi_\lambda(H_\alpha)=c\lambda(H_\alpha)\neq 0,$
    which implies the assertion of the lemma.

    Let us now show that the inverse to Lemma 11 is also true.

    Let  $\lambda\in a^*$ be a dominant weight, $\pi_\lambda$
    the representation of ${\cal U}(\underline{g})$ in the module
    $E^\lambda$ with the highest weight $\lambda$ and the highest weight
    vector $v$ of the length 1.

{\bf Lemma 12.} {\it Let $\alpha\in\triangle^+(\lambda)$, then
$\pi_\lambda(X_{-\alpha})v=0$.}

{\it Proof.} Consider $E^\lambda$ as a module over the $sl_2$-subalgebra
of $\underline{g}$, generated by the elements $X_{-\alpha},H_\alpha,X_\alpha$.
Assume $\pi_\lambda(X_{-\alpha})v\neq 0$.
Since $\pi_\lambda(X_\alpha X_{-\alpha})v=
\pi_\lambda(X_{-\alpha} X_\alpha)v+\pi_\lambda(H_\alpha)v=
\lambda(H_\alpha)v=0$, and\\
 $\pi_\lambda(H_\alpha X_{-\alpha})v=
\pi_\lambda(X_{-\alpha} H_\alpha)v+\pi_\lambda(-2X_{-\alpha})v=
-2\pi_\lambda(X_{-\alpha})v$, then $\pi_\lambda(X_{-\alpha})v$
generates an infinite dimensional $sl_2$-submodule of $E^\lambda$,
which  contradicts the fact that $E^\lambda$
is finite dimensional.

{\bf Lemma 13.} {\it The highest weight vector $v\in E^\lambda$
is an eigenvector for the subgroup $K^\lambda\subset K$.}

Since $K^\lambda$ is the integral subgroup of the subalgebra
$k^\lambda\subset k$, it is enough to check the assertion of the lemma
for the subalgebra $k^\lambda$. It is generated by the subalgebra $m$,
for which $v$ is an eigenvector, and by the elements
$X_\alpha-X_{-\alpha}$ and $i(X_\alpha+X_{-\alpha}),\
\alpha\in\triangle^+(\lambda)$, which annihilate $v$ according to Lemma 12.

{\bf Lemma 14.} {\it Let $f(k)$ be the covariant symbol of an operator
$A\in End\ E^\lambda$. Then the function $f(k)$ is invariant under the
right shifts by the elements of $K^\lambda$.}

{\it Proof.} It follows from Lemma 13 and a unitarity of
the representation of $K$ in $E^\lambda$, that for
$l\in K^\lambda\quad lv=cv$, where $c$ is a constant of module 1.
Therefore, $f(kl)=<Aklv,klv>=<Akv,kv>=f(k)$. The lemma is proved.

{\bf Theorem 2.} {\it All the functions from
$A(\lambda+\rho),\ \lambda\in a^*,$ are invariant under the
right shifts by the elements of $K^\lambda$ iff $\lambda$
is relatively regular.}

{\it Proof.} The necessity is proved in Lemma 11. Let us prove
the sufficiency. Assume that $\lambda$ is regular relative to
$\sigma\in\Sigma$. Denote by $a^*_\sigma$ the set of $\mu\in a^*$
such that
$\mu(H_\alpha)=0$ for $\alpha\in\sigma$. In particular,
$\lambda\in a^*_\sigma$. Denote by $\omega^\alpha$
the fundamental weight, corresponding to $\alpha\in\Sigma$
(i.e. $\omega^\alpha(H_\beta)=\delta^\alpha_\beta$ for
$\alpha,\beta\in\Sigma$, where $\delta^\alpha_\beta$ is the
Kronecker symbol). The fundamental weights
$\omega^\alpha,\ \alpha\in\Sigma\backslash\sigma$, form a base in
$a^*_\sigma$, and their linear combinations with coefficients in {\bf N}
are the dominant weights, which are regular relative to $\sigma$.
Consider $s_\mu(k)$ and $s_\mu(kl)$ for fixed $k\in K,l\in K^\lambda$
and $u\in{\cal U}(\underline{g})$ as polynomials in $\mu$ on $a^*_\sigma$.
It follows from Proposition 2 and Lemma 14, that these polynomials
coincide on the dominant weights $\mu$, which are regular relative to
$\sigma$. It is easy to conclude from the above, that these polynomials
coincide identically, which proves the theorem.

{\bf Corollary.} {\it Let $\Omega_\lambda\subset k$ be an adjoint orbit
of the point $iH_\lambda$. The algebra $A(\lambda+\rho)$ can be pushed
forward to $\Omega_\lambda$ by the mapping $\Psi_\lambda$
iff $\lambda$ is relatively regular.}

{\bf 7. Quantization on flag manifolds}

Let $d$ be a nonnegative integer. Denote by
${\cal U}_d(\underline{g})$ the subspace of ${\cal U}(\underline{g})$,
generated by all monomials of the form $X_1\dots X_j$, where
$X_1,\dots,X_j\in\underline{g}$ and $j\leq d$. The subspaces
$\{{\cal U}_d(\underline{g})\}$ determine the canonical filtration on
${\cal U}(\underline{g})$ (see \cite{Dix}).

The symmetric algebra $S(\underline{g})$ can be identified with the
space of polynomials on $k$, so that the element $X\in\underline{g}$
corresponds to the linear functional on $k$,
$\tilde X(Y)=(X,Y),\ Y\in k$. Let $S^d(\underline{g})$
be the space of homogeneous polynomials on $k$ of degree $d$.
The graded algebra, associated with the canonical filtration
is canonically isomorphic to $S(\underline{g})$, so that
${\cal U}_d(\underline{g})/{\cal U}_{d-1}(\underline{g})$
corresponds to $S^d(\underline{g})$. For $u\in{\cal U}_d(\underline{g})$
let $\underline{u}^{(d)}$ denote the corresponding element of
$S^d(\underline{g})$.
Usually we will omit the superscript $d$ in $\underline{u}^{(d)}$.
If $u=X_1\dots X_d\in{\cal U}_d(\underline{g})$, then
$\underline{u}=\tilde X_1\dots\tilde X_d$.
We say, that a monomial $u=X_1\dots X_d\in{\cal U}(\underline{g})$
has a canonical form, if for some integers $k,l$, such that
$0\leq k\leq l\leq d$, an element $X_j$ belongs to
$\underline{n}^-,\underline{h},\underline{n}$ if
$0<j\leq k,\ k<j\leq l,\ l<j\leq d$, respectively.
Recall, that the projection of $u\in{\cal U}(\underline{g})$
onto ${\cal U}(\underline{h})$
in (2) is denoted by $u_0$.

{\bf Proposition 3.} {\it Let $t\in {\bf R},\ \lambda\in a^*,\
u\in{\cal U}_d(\underline{g}),\ k\in K$.
Then }
$$
\lim_{t\to\infty}t^{-d}s_{t\lambda}u(k)=
i^{-d}\underline{u}^{(d)}(\Psi_\lambda(k)).
$$

{\it Proof.} Since $s_\lambda,\ \Psi_\lambda$ and the mapping
$u\mapsto\underline{u}$ are $K$-equivariant, it is enough to check
the equality in the proposition only at the point $e\in K$,
i.e. one has to show, that
\begin{equation}
\lim_{t\to\infty}t^{-d}u_0(t\lambda)=
i^{-d}\underline{u}(iH^\lambda).
\end{equation}

Assume that $u\in{\cal U}_{d-1}(\underline{g})$.
Then $u_0(t\lambda)$ is a polynomial in $t$ of degree less than $d$,
hence the left-hand side in (7) is 0. On the other hand
$\underline{u}=0$, since ${\cal U}_{d-1}(\underline{g})$
is the kernel of the mapping $u\mapsto\underline{u}^{(d)}$.
Now it remains to check (7) for a monomial
$u=X_1\dots X_d\in{\cal U}_d(\underline{g})$, which has
a canonical form. Since
\begin{equation}
i^{-d}\underline{u}(iH^\lambda)=
i^{-d}\tilde X_1(iH^\lambda)\dots\tilde X_d(iH^\lambda)=
(X_1,H^\lambda)\dots(X_d,H^\lambda),
\end{equation}
the expression in (8) for $u\in{\cal U}(\underline{h})$ equals to
$u_0(\lambda)$ . In this case $u_0(\lambda)$ is homogeneous of degree $d$,
so the expression under the limit on the left-hand side of (7)
does not depend on $t$ and the left-hand side of (7) is also
equal to $u_0(\lambda)$. Assume now that $u\notin{\cal U}(\underline{h})$.
Then the factor $(X_j,H^\lambda)$ on the right-hand side of (8), for which
$X_j\notin\underline{h}$, is equal to 0. On the other hand, in this case
$u_0=0$, which concludes the proof.

Let $V$ be a $K$-module, $\nu\in a^*$ a dominant weight. Denote by $V^\nu$
the isotropic component of $V$, corresponding to the representation
with the highest weight $\nu$.
Let $\Omega\subset k$ be an adjoint orbit of the group $K$.
Denote by $Q(\Omega)$ the space of $K$-finite functions on $\Omega$,
and by $R(\Omega)$ the space of regular functions on $\Omega$,
i.e., the space of restrictions of the functions from
$S(\underline{g})$ on $\Omega$.

{\bf Lemma 15.} {\it (i) The spaces $Q(\Omega)$ and $R(\Omega)$ coincide. \\
(ii)  The space $R(\Omega)^\nu$ for each dominant weight $\nu$
is finite dimensional. \\
(iii) There exist homogeneous polynomials
$\psi_1,\dots,\psi_n\in S(\underline{g})^\nu$ such that their
restrictions on $\Omega$ form a base in $R(\Omega)^\nu$.}

{\it Proof.} Since each element in $S(\underline{g})$ is $K$-finite,
 $R(\Omega)\subset Q(\Omega)$. It follows from the Frobenius reciprocity
theorem, that  the irreducible representation of $K$
with the highest weight $\nu$  has a finite multiplicity in $Q(\Omega)$.
Hence $Q(\Omega)^\nu$ and, therefore, $R(\Omega)^\nu$ are finite dimensional.
Assume that there exists a function
$f\in Q(\Omega)^\nu\backslash R(\Omega)^\nu$ for some $\nu$.
Then $f$ is orthogonal to $R(\Omega)$ in $L^2(\Omega, dx)$, where
$dx$ is the $K$-invariant measure on $\Omega$.
It follows from the Stone-Weierstrass theorem, that $R(\Omega)$
is dense in $C(\Omega)$ and, therefore, in $L^2(\Omega,dx)$,
thus $f=0$, which implies (i). Finally, since
$S(\underline{g})^\nu=\oplus (S^d(\underline{g}))^\nu$, the polynomials
from (iii) can be constructed inductively.

{\bf Lemma 16.} {\it Let $\lambda\in a^*$ be a relatively regular element,
$\nu\in a^*$ a dominant weight,
 $u\in{\cal U}_d(\underline{g})^\nu,\ t\in{\bf R},\
\{f_1,\dots,f_n\}$ a base in $R(\Omega_\lambda)^\nu$.
Then $s_{t\lambda}u$ can be realized as a linear combination of functions
$f_j\circ\Psi_\lambda,\ j=1,\dots,n$,
with  coefficients, which are polynomial in $t$
of a degree not greater than $d$.}

{\it Proof.} It is easy to notice that for a given $\lambda$
 the stabilizer $K^{t\lambda}$ coincides with
$K^{\lambda}$ for all real $t\neq 0$, and with $K$ for $t=0$.
Since $t\lambda$ is relatively regular for all $t$, it follows from
Theorem 2 that $s_{t\lambda}u$ is the pullback of some function
$\psi_t$ on $\Omega_\lambda$, that is,
$s_{t\lambda}u=\psi_t\circ\Psi_\lambda$.
The $K$-equivariance of $s_\lambda$ and Lemma 15 imply that
$\psi_t\in R(\Omega_\lambda)^\nu$, thus there exist functions
$a_j(t),\ j=1,\dots,n,$ such that $\psi_t=\Sigma a_j(t)f_j$.
Since $R(\Omega_\lambda)^\nu$ is $n$-dimensional,
there exist points
$x_1,\dots,x_n\in\Omega_\lambda$ such that the matrix $(f_i(x_j))$
is invertible. Let $(b_{ij})$ be its inverse, then
$a_i(t)=\Sigma b_{ij}\psi_t(x_j)$. Now the assertion of the lemma
follows from the fact that for each $k\in K$ the function
$\psi_t(\Psi_\lambda(k))=s_{t\lambda}u(k)$  is polynomial
in $t$ of a degree not greater than $d$.

Let $\beta:S(\underline{g})\to{\cal U}(\underline{g})$
be the symmetrization mapping (see \cite{Dix}). It is $K$-equivariant,
and for a homogeneous polynomial $\psi\in S^d(\underline{g})$
holds $\beta \psi\in{\cal U}_d(\underline{g})$
and $\underline{\beta \psi}=\psi$.

{\bf Lemma 17.} {\it Let $\lambda\in a^*$ be relatively regular,
 $f\in R(\Omega_\lambda)$. There exist
 $u_j\in{\cal U}_{d(j)}(\underline{g}),\ j=1,\dots,m,$
 and rational functions $\{a_j(t)\}$ with no poles at infinity,
 such that for all $t\neq 0$ different from the poles of the functions
 $\{a_j(t)\}$,  holds}
$$
f\circ\Psi_\lambda=\sum a_j(t)t^{-d(j)}s_{t\lambda}u_j.
$$

{\it Proof.} It is enough to state the lemma for functions
$f_1,\dots,f_n$ which\\
 form the base of
$R(\Omega_\lambda)^\nu$ for some $\nu$.
According to Lemma 15, one can assume that $f_j$ is a restriction
of a homogeneous polynomial  $\psi_j$ of degree $d(j)$.
It follows from Lemma 16 that there exist polynomials
 $b_{jl}(t)$ of degree not greater than $d(j)$, such that
\begin{equation}
    (\frac{i}{t})^{d(j)}s_{t\lambda}\beta\psi_j=
      \sum_l t^{-d(j)}b_{jl}(t)f_l\circ\Psi_\lambda.
\end{equation}
Using Proposition 3 and the properties of $\beta$, pass to the limit
 on both sides of (9) as $t\to\infty$. Since the left-hand
side of (9) tends to $f_j\circ\Psi_\lambda$, then the matrix
of rational functions $(t^{-d(j)}b_{jl}(t))$ tends to the identity matrix as
$t\to\infty$. Therefore its determinant is not identically zero,
and there exists its inverse matrix of rational functions,
$(a_{lj}(t))$, which also tends to the identity matrix as
$t\to\infty$. The assertion of the lemma is obtained by applying
the matrix $(a_{lj}(t))$ to both sides of (9) and setting
$u_j=i^{d(j)}\beta\psi_j$.

It is known that there exists a natural Poisson structure
on the symmetric algebra $S(\underline{g})$, which is considered
as the algebra of polynomials on $k$ . It is defined by the
(real) Poisson bracket $\{\cdot,\cdot\}$ on $k$ such that for
$X,Y\in\underline{g},\ Z=[X,Y],$ holds
$\{\tilde X,\tilde Y\}=\tilde Z$. The symplectic leaves of that
Poisson structure are the (co)adjoint orbits of $K$ in $k$
(the Lie algebra $k$ is identified with its dual space
$k^*$ via the Killing form). The symplectic structure on a
(co)adjoint orbit $\Omega\subset k$ is defined by the Kirillov
symplectic form, and the corresponding Poisson bracket
$\{\cdot,\cdot\}_\Omega$ is the restriction of
$\{\cdot,\cdot\}$ on $\Omega$, i.e. for
$f,g\in S(\underline{g})$ holds $\{f,g\}|_\Omega
                           =\{f|_\Omega,g|_\Omega\}_\Omega$ (see \cite{BF}).

The multiplication and the Poisson bracket on  $S(\underline{g})$
are connected with the canonical filtration on
${\cal U}(\underline{g})$ in the following way (see \cite{Dix}). Let
$u_\varepsilon\in{\cal U}_{d_\varepsilon}(\underline{g}),\
 \varepsilon=1,2,\ d=d_1+d_2$, then $u_1u_2\in{\cal U}_d(\underline{g}),\
u_1u_2-u_2u_1\in{\cal U}_{d-1}(\underline{g})$,  and
\begin{equation}
(\underline{u_1u_2})^{(d)}=\underline{u_1}^{(d_1)}
\cdot\underline{u_2}^{(d_2)},\quad
(\underline{u_1u_2-u_2u_1})^{(d-1)}=
\{\underline{u_1}^{(d_1)},\underline{u_2}^{(d_2)}\}.
\end{equation}

Let $\lambda\in a^*$ be relatively regular,
$t\in {\bf R}\backslash\{0\}$.
It follows from  the fact that $K^{t\lambda}=K^\lambda$,
and Theorem 2, that the algebra $A(t\lambda+\rho)$ can be pushed forward
to the orbit $\Omega_\lambda$ via $\Psi_\lambda$. Denote by
${\cal A}_{1/t}^{(\lambda)}$ and $*_{1/t}^{(\lambda)}$
the pushforward of $A(t\lambda+\rho)$ on $\Omega_\lambda$,
and the multiplication in it, respectively.
Usually the superscript $\lambda$ in these notations will be omitted.
Since each element in $A(t\lambda+\rho)$ is $K$-finite,
Lemma 15 implies that
${\cal A}_{1/t}^{(\lambda)}\subset R(\Omega_\lambda)$.

{\bf Theorem 3.} {\it Let $\lambda\in a^*$ be relatively regular,
$t\in {\bf R}\backslash\{0\}$.\\
(i) Each function $f\in R(\Omega_\lambda)$ belongs to
${\cal A}_{1/t}^{(\lambda)}$ for
all but a finite number of values of $t$.\\
(ii) For $f_1,f_2\in R(\Omega_\lambda),\
 x\in\Omega_\lambda,$ the product
$(f_1*_{1/t}f_2)(x)$ is a rational function of $t$,
with no pole at infinity.\\
(iii) For $f_1,f_2\in R(\Omega_\lambda)$
 the limits below hold}
$$
\lim_{t\to\infty} f_1*_{1/t}f_2=f_1f_2;\quad
\lim_{t\to\infty} t(f_1*_{1/t}f_2-f_2*_{1/t}f_1)=
      i\{f_1,f_2\}_{\Omega_\lambda}.
$$

{\it Proof.} Assertion (i) directly follows from Lemma 17.
In accordance with Lemma 17, there exist the elements
$u_{\varepsilon j}\in{\cal U}_{d_\varepsilon(j)}(\underline{g}),\
\varepsilon=1,2,\ j=1,\dots,n(\varepsilon),$
and the rational functions $a_{\varepsilon j}(t)$,
with no pole at infinity, such that for all $t\neq 0$, different
from the poles of $\{a_{\varepsilon j}(t)\}$, holds
\begin{equation}
              f_\varepsilon\circ\Psi_\lambda=\sum_{j=1}^{n(\varepsilon)}
              a_{\varepsilon j}(t) t^{-d_\varepsilon(j)}
              s_{t\lambda}u_{\varepsilon j}.
\end{equation}

Passing to the limit in (11) as $t\to\infty$, and denoting
$\lim_{t\to\infty}a_{\varepsilon j}(t)=a_{\varepsilon j}$, one gets from
Proposition 3, that
\begin{equation}
              f_\varepsilon\circ\Psi_\lambda=\sum_j
              a_{\varepsilon j}i^{-d_\varepsilon(j)}
              \underline{u}_{\varepsilon j}\circ\Psi_\lambda.
\end{equation}

 Let us calculate the product  $f_1*_{1/t}f_2$,
pulling it back to $K$ via $\Psi_\lambda$:
\begin{equation}
                   (f_1*_{1/t}f_2)\circ\Psi_\lambda
              =\sum_{j,l}a_{1j}(t)a_{2l}(t)t^{-d(j,l)}
                     s_{t\lambda }(u_{1j}u_{2l}),
\end{equation}
where $d(j,l)=d_1(j)+d_2(l)$.
Since $u_{1j}u_{2l}\in{\cal U}_{d(j,l)}(\underline{g})$, for each
$k\in K$\\
$s_{t\lambda }(u_{1j}u_{2l})(k)$
is a polynomial in $t$ of a degree not greater than $d(j,l)$.
Thus each summand in (13) is a rational function of $t$,
with no pole at infinity, which proves (ii).

Pass to the limit  in (13) as $t\to\infty$, using Proposition 3
and taking into account (10), (11) and (12):
$$
\lim_{t\to\infty}(f_1*_{1/t}f_2)\circ\Psi_\lambda=
              \sum_{j,l}a_{1j}a_{2l}i^{-d(j,l)}
              (\underline{u_{1j}u_{2l}})\circ\Psi_\lambda=
$$
$$
             ( \sum_{j,l}a_{1j}a_{2l}
        i^{-d_1(j)}\underline{u}_{1j}i^{-d_2(l)}\underline{u}_{2l})
        \circ\Psi_\lambda=
              (f_1f_2)\circ\Psi_\lambda.
$$
Thus we obtain the first of the two limits in (iii).
The second limit is calculated similarly:
$$
\lim_{t\to\infty}t(f_1*_{1/t}f_2-f_2*_{1/t}f_1)\circ\Psi_\lambda=
\lim_{t\to\infty}\sum_{j,l}a_{1j}(t)a_{2l}(t)t^{-d(j,l)+1}
                     s_{t\lambda }(u_{1j}u_{2l}-u_{2l}u_{1j})=
$$
$$
 \sum_{j,l}a_{1j}a_{2l}i^{-d(j,l)+1}
(\underline{u_{1j}u_{2l}-u_{2l}u_{1j}})\circ\Psi_\lambda=
$$
$$
 i(\sum_{j,l}a_{1j}a_{2l}\{
        i^{-d_1(j)}\underline{u}_{1j},i^{-d_2(l)}\underline{u}_{2l}
                \})\circ\Psi_\lambda=
    i\{f_1,f_2\}_{\Omega_\lambda}\circ\Psi_\lambda.
$$

{\it Remark 1.} The assertion (i) of Theorem 3 means, that the algebras
${\cal A}_{1/t}^{(\lambda)}$ have the weak nesting property.
It is proved in \cite{RCG} that if $\lambda$ is a dominant weight
and  $n\in{\bf N},\
{\cal A}_{1/n}^{(\lambda)}\subset{\cal A}_{1/(n+1)}^{(\lambda)}$,
i.e., the algebras of covariant symbols,
$\{{\cal A}_{1/n}^{(\lambda)}\}$, have the strong nesting property.

{\it Remark 2.} If $\lambda$ is a dominant weight,
the assertion (ii) of Theorem 3 is proved
when $\Omega_\lambda$ is a compact hermitian symmetric space,
and conjectured for an arbitrary generalized flag manifold in \cite{CGR}.

{\it Remark 3.} The assertion (iii) of Theorem 3
is proved in \cite{CGR} for the case when $\lambda$
is a dominant weight.

{\it Remark 4.} The connection between the $K$-orbit of the
highest weight vector $v\in E^\lambda$
and the adjoint orbit $\Omega_\lambda$ with their application to
covariant symbols is studied in \cite{Wild}.

Now we can define the Berezin's quantization
on a generalized flag manifold, using the algebras of mixed symbols.
Let $\lambda\in a^*$ be a relatively regular weight,
and there exists an element $w\in W$ such that the weights
$w\lambda$ and $w\cdot\lambda$ are dominant. Then for all natural
$n$ the weights $w\cdot(n\lambda)=(n-1)w\lambda+w\cdot\lambda$
are dominant.
The algebra $A(n\lambda+\rho)$ has the representation $\tau_{n\lambda}$
in the $\underline{g}$-module $E^{w\cdot(n\lambda)}$. Pushing forward
$A(n\lambda+\rho)$ together with its representation $\tau_{n\lambda}$
to $\Omega_\lambda$, we obtain the algebra
${\cal A}_{1/n}^{(\lambda)}$ with a representation in
$E^{w\cdot(n\lambda)}$. Denote
$\Omega=\Omega_\lambda,\ F=\{1,1/2,1/3,\dots\},\ h=1/n\in F,\ {\cal A}_h=
{\cal A}_{1/n}^{(\lambda)},\ H_h=E^{w\cdot(n\lambda)}$.
Then Theorem 3 implies

{\bf Theorem 4.} {\it The algebras $\{{\cal A}_h\},\ h\in F$,
together with their representations in $H_h$
define the Berezin's quantization on the (co)adjoint orbit
$\Omega$ of the group $K$, which is a generalized flag manifold.}


\begin{thebibliography}{99}
\bibitem{Ber1} F.A.Berezin, {\it Quantization on bounded complex
domains}, Soviet Math. Dokl., {\bf 211}(1973), 1263-1266.
\bibitem{Ber2} F.A.Berezin, {\it}  {\it General concept of quantization},
Comm. Math. Phys. {\bf 40} (1975), 153-174.
\bibitem{CGR} M.Cahen, S.Gutt, J.Rawnsley, {\it Quantization of
K\"ahler manifolds II}, TAMS, {\bf 337} (1993), 73-98.
\bibitem{Wild} N.Wildberger, {\it On the Fourier transform of a
compact semisimple Lie group}, Preprint of the Univ. of N.S.W., Australia.
\bibitem{Duf}   M.Duflo,
{\it Representations irreductibles des groupes semi-simples complexes},
 Lecture Notes in Math., vol. 497, Springer-Verlag, Berlin,
Heidelberg and New York, 1975.
\bibitem{Dix} J.Dixmier, {\it Alg\`ebres enveloppantes}
Gauthier-Villars, Paris, 1974.
\bibitem{BF} F.Bayen et al.,{\it Deformation theory and quantization},
Ann. Phys. {\bf 111} (1978) 61-151.
\bibitem{RCG}   J.Rawnsley, M.Cahen, S.Gutt, {\it Quantization of K\"ahler
manifolds I:  geometric interpretation of Berezin's quantization},
 J. Geom. Phys. {\bf 7} (1990),  45-62.
\end{thebibliography}
\end{document}